\documentclass[reqno,a4paper,11pt]{amsart}

\usepackage[utf8]{inputenc}
\usepackage{graphicx}
\usepackage{subfigure}
\usepackage{amsmath, amssymb, amsthm, mathtools, braket,times}
\usepackage{amsaddr}
\usepackage{enumitem}
\usepackage{cite}
\usepackage{graphicx}
\usepackage{array}
\usepackage{tikz}
\usepackage{mathrsfs}
\usetikzlibrary{arrows}
\usetikzlibrary{shapes}
\usepackage{cleveref}
\usepackage[english]{babel}
\usepackage[T1]{fontenc}
\usepackage{fourier}
\usepackage{bbm}
\usepackage{color}
\usepackage{fullpage}

\newtheorem{thm}{Theorem}[section]

\newtheorem{lem}[thm]{Lemma}

\newtheorem{conj}[thm]{Conjecture}

\newtheorem{problem}[thm]{Problem}

\numberwithin{equation}{section}

\newcommand{\D}{\mathrm{deg}}

\title{Equitable vertex arboricity conjecture holds for graphs with low degeneracy}
\author{Xin Zhang, ~~Bei Niu, ~~Yan Li, ~~Bi Li}
\address{School of Mathematics and Statistics, Xidian University, Xi'an, 710071, China}
\email{xzhang@xidian.edu.cn (X.\,Zhang)}
\thanks{National Natural Science Foundation of China (Nos.\,11871055,11701440)}

\begin{document}

\begin{abstract}
The equitable tree-coloring can formulate a structure decomposition problem on the communication network with some security considerations. Namely,
an equitable tree-$k$-coloring of a graph is a vertex coloring using $k$ distinct colors such that every color class induces a forest and the sizes of any two color classes differ by at most one.
In this paper, we show some theoretical results on the equitable tree-coloring of graphs by proving that
every $d$-degenerate graph with maximum degree at most $\Delta$ is equitably tree-$k$-colorable for every integer $k\geq (\Delta+1)/2$ provided that $\Delta\geq 9.818d$, confirming the equitable vertex arboricity conjecture for graphs with low degeneracy.
\end{abstract}

\maketitle

\section{Introduction}

Almost all relationships in the real world can be modeled by networks, so doing theoretical research on the structure of a network is necessary and interesting. Actually, finding treelike structures in a network is a fascinating topic. It attracts the attentions from many researchers including \cite{AD16,HHZ.2018,L18,L19,MCX16}.

While decomposing a large communication network into small pieces, we are sometimes required to make each of its small pieces having an ``acyclic'' property due to some security considerations. In this way,
the local structure around a node in a small piece is so clear that it can be easily tested using some classic algorithmic tools, and thus we can effectively identify any possible node failure.
Meanwhile, we do not want to have too many small pieces and hope that the scales of any two distinct pieces would not differ a lot. These requirements help us to keep up the communication network effectively even if we have a low budget.
Practically, this structure decomposition problem can be modeled by \emph{equitable tree-coloring}, which was introduced by Wu, Zhang and Li \cite{WZL.2013} in 2013.

From now on, we do not distinguish graph from network, and introduce some graph-based definitions and notations. A graph is \emph{$d$-degenerate} if its every subgraph has a vertex of degree at most $d$. The smallest integer such that $G$ is $d$-degenerate is the \emph{degeneracy} of $G$, which is used to measure the sparseness of the graph. Clearly, the degeneracy of a graph is upper-bounded by its maximum degree. Let $V_1$ and $V_2$ be two subsets of the vertex set $V(G)$ of a graph $G$. By $e(V_1,V_2)$, we denote the number of edges that have one end-vertex in $V_1$ and the other in $V_2$. By $\D_G(v)$, we denote the degree of the vertex $v$ in $G$. By $e(G)$, we denote the number of edges of $G$. For other undefined notation, we refer the readers to the book \cite{Bondy.2008}.

An \emph{equitable tree-$k$-coloring} (resp.\,\emph{equitable $k$-coloring}) of $G$ is a function $c$ from $V(G)$ to $\{1,2,\ldots,k\}$ so that $c^{-1}(i)$ induces a forest (resp.\,an independent set) for every $1\leq i\leq k$, and $\big||c^{-1}(i)|-|c^{-1}(j)|\big|\leq 1$ for every pair of $i,j$ with $1\leq i<j\leq k$. The smallest integer $k$ such that $G$ admits an equitable tree-$k$-coloring, denoted by $va_{eq}(G)$, is the \emph{equitable vertex arboricity} of $G$. An equitably tree-$k$-colorable graph may not be equitably tree-$k'$-colorable for some $k'>k$.
For instance, $va_{eq}(K_{9,9})=2$, but $K_{9,9}$ is not equitably tree-3-colorable.
In view of this, a new chromatic parameter, \emph{the  equitable vertex arboreal threshold},  was introduced. Precisely, it is the smallest integer $k$, denoted by $va_{eq}^*(G)$, such that $G$ has an equitable tree-$k'$-coloring for any integer $k'\geq k$.
Those concepts were initially introduced by Wu, Zhang and Li \cite{WZL.2013} in 2013, who proposed the following two conjectures.

\begin{conj}\label{conj1}[Equitable Vertex Arboricity Conjecture (EVAC)]
Every graph with maximum degree at most $\Delta$ is equitably tree-$k$-colorable for any integer $k\geq (\Delta+1)/2$, i.e, $va_{eq}^*(G)\leq \lceil(\Delta+1)/2\rceil$.
\end{conj}

\begin{conj}\label{conj2}
There exists a constant $k$ such that every planar graph is equitably tree-$k'$-colorable for any integer $k'\geq k$, i.e, $va_{eq}^*(G)\leq k$.
\end{conj}

Zhang \cite{Z.2016} showed that every subcubic graph is equitably tree-$k$-colorable for any integer $k\geq 2$, and thus Conjecture \ref{conj1} holds for $\Delta\leq 3$.
Chen \emph{et al\,} \cite{CGSWW.2017} verified Conjecture \ref{conj1} for all 5-degenerate graphs, which implies that Conjecture \ref{conj1} holds for $\Delta\leq 5$. On the other hand, Zhang and Wu \cite{ZW.2014}, and Zhang and Niu \cite{ZN.2020} proved that Conjecture \ref{conj1} holds for $\Delta\geq (|V(G)|-1)/2$.

In 2015, Esperet, Lemoine and Maffray \cite{ELM.2015} confirmed Conjecture \ref{conj2} by showing that
$va_{eq}^*(G)\leq 4$ for any planar graph $G$. Since the vertex arboricity of every planar graph is at most 3 \cite{CK}, it is natural to put forward the following
conjecture.

\begin{conj}\label{conj3}
Every planar graph is equitably tree-$3$-colorable, i.e., $va_{eq}^*(G)\leq 3$ for any planar graph $G$.
\end{conj}

Conjecture \ref{conj3} has been verified for some classes of graphs including
planar graphs with maximum degree at most 5 \cite{CGSWW.2017},
planar graphs with girth at least 5 \cite{WZL.2013}, and planar graphs with some conditions on the cycles \cite{Z.2015}.

For $d$-generate graphs $G$, Esperet, Lemoine and Maffray \cite{ELM.2015} showed that $va_{eq}^*(G)\leq 3^{d-1}$.
This implies the following result.
\begin{thm}\label{thmold}
 Every $d$-degenerate graph with maximum degree at most $\Delta$ is equitably tree-$k$-colorable for any integer $k\geq (\Delta+1)/2$ provided that $\Delta\geq 2 \cdot 3^{d-1}-1$.
\end{thm}

The aim of this paper is to lower the exponential lower bound of $\Delta$ in Theorem \ref{thmold} (for $d\geq 4$) to $9.818d$. Precisely, we are to prove the following theorem, confirming the equitable vertex arboricity conjecture for graphs with low degeneracy in terms of its maximum degree

\begin{thm}\label{main}
Every $d$-degenerate graph with maximum degree at most $\Delta$ is equitably tree-$k$-colorable for any integer $k\geq (\Delta+1)/2$ provided that $\Delta\geq 9.818d$.
\end{thm}

\section{The proof of Theorem \ref{main}}

To begin with, it may be worthy to mention a relative result as follows.

\begin{thm}\cite[Kostochka and Nakprasit]{KN.2003}\label{thmKN}
\textit{Every $d$-degenerate graph with maximum degree at most $\Delta$ is equitably $k$-colorable for any integer $k\geq (\Delta+d+1)/2$ provided that $\Delta\geq 27d$ and $d\geq 2$.}
\end{thm}

Remark that Theorem \ref{main} cannot be established by Theorem \ref{thmKN} or by its proof with just replacing the condition on $\Delta$.
Although one may see below that the idea of proving Theorem \ref{main} partially comes from the proof of Theorem \ref{thmKN}, one shall simultaneously realize, on the other hand, that some quite different ideas are introduced to break a barrier. 

Actually, an equitable tree-$k$-coloring of a graph $G$ on $n$ vertices with $k\nmid n$, say $n=kt-r$ where $1\leq r\leq k-1$, cannot be directly reduced from an equitable tree-$k$-coloring of the graph $G\cup K_r$ (the reason is that in an tree-$k$-coloring of $K_r$, we do not have the good property that distinct vertices receive different colors --- but this good property exists and is applicable, as in the proof of Theorem \ref{thmKN}, if one is working on the equitable $k$-coloring).
Therefore, in our proof, if we still want to reduce the equitable tree-$k$-coloring of $G\cup K_r$ to $G$, we shall fix in advance a rainbow coloring for that $K_r$, at the expense of taking care of the parameter $r$ throughout the remaining arguments, see \eqref{bound-r}, \eqref{i-1}, \eqref{proc} and \eqref{basic} for some undertanding later.

In what follows, we give the detailed proof of Theorem \ref{main}. During the proof, the following lemma is  applicable for setting up.

\begin{lem}\cite[Kostochka and Nakprasit]{KN.2003}
\label{lem-1}
Let $G$ be a $d$-degenerate graph on $n$ vertices, where $d\geq 2$ is an integer. If
$v_1,v_2,\ldots,v_n$ is an ordering of  the vertices of $G$ from the highest degree to the lowest degree, then $\D_G(v_i)<d(1+\frac{n}{i})$ for every $i=1,2,\ldots,n$.
\end{lem}

\proof[\emph{\textbf{Proof of Theorem \ref{main}}}] Since 0-degenerate graphs or 1-degenerate graphs are forests and it is trivial that every forest admits an equitably tree-$k$-colorable for every integer $k\geq 1$, we assume $d\geq 2$.

We proceed by induction on $|V(G)|$. If $|V(G)|\leq 2$, then the result is trivial.
So we assume that the result holds for every $d$-degenerate graph on less than $n$ vertices, and work on a $d$-degenerate graph on $n$ vertices with $k(t-1)< n\leq kt$, where $t$ is a positive integer.

If $t=2$ (resp.\,$t=1$), then we arbitrarily partition the vertices of $G$ into $k$ subsets such that each of them has either one or two vertices (resp.\,at most one vertex). Clearly, this partition is an equitable tree-$k$-coloring of $G$. Hence we assume $t\geq 3$, and precisely, let $n=kt-r$, where $0\leq r< k$. We now claim that
\begin{align}\label{bound-r}
2(r+1)\leq d.
\end{align}

Since $G$ is $d$-degenerate, there is a vertex $x$ with $\D_G(x)\leq d$. By the induction hypothesis, $G-x$ has an
equitable tree-$k$-coloring $c$. Since $|V(G-x)|=kt-(r+1)\geq k(t-1)$, exactly $r+1$ color classes of $c$ are of size $t-1$ and other color classes are of size $t$. If $d<2(r+1)$, then among those $r+1$ color classes of size $t-1$, there is at least one having at most one neighbor of $x$. In such a case we move $x$ into this color class and get an equitably tree-$k$-coloring of $G$. Hence we assume \eqref{bound-r}.

Suppose that
$v_1,v_2,\ldots,v_n$ is an ordering of the vertices of $G$ from the highest degree to the lowest degree
with the property that
\begin{align}\label{ordering}
    \D_G(v_{1})\geq \cdots \geq \D_G(v_{\mu})\geq \lambda(t) d>\D_G(v_{\mu+1})\geq \cdots \geq \D_G(v_{n}),
\end{align}
where 
\begin{align}\label{lambda}
\lambda:=\lambda(t)=1+\frac{t}{s(t)}
\end{align}
and
\begin{align}\label{s(t)}
s:=s(t)=\left\{
          \begin{array}{ll}
            1, & \hbox{if $3\leq t\leq 6$;} \\
            \big\lceil\frac{t}{5}\big\rceil, & \hbox{if $t\geq 7$.}
          \end{array}
        \right.
\end{align}
Note that
\begin{align}
\label{lmdasmall} 4\leq \lambda&=1+t\leq 7~{\rm if}~3\leq t \leq 6,~{\rm and}\\
\label{lmda} 4.1818<6-\frac{20}{7+4} \leq 6-\frac{20}{t+4}=1+\frac{t}{\frac{t+4}{5}}\leq \lambda&=1+\frac{t}{\big\lceil\frac{t}{5}\big\rceil}\leq 1+\frac{t}{\frac{t}{5}}=6~{\rm if}~t\geq 7.
\end{align}

Since $\D_G(v_\mu)<d(1+\frac{n}{\mu})$ by Lemma \ref{lem-1}, 
$\mu< \frac{n}{\lambda-1}< n$ by \eqref{ordering}, \eqref{lmdasmall} and \eqref{lmda}.
This implies, by the induction hypothesis, that the graph induced by $V'=\{v_1,v_2,\ldots,v_{\mu}\}$ has an equitable tree-$k$-coloring $c'$, and
\begin{align}\label{s}
  \textit{every~color~class~of}~c'~\textit{contains~at~most}~\bigg\lceil\frac{\mu}{k}\bigg\rceil\leq \bigg\lceil\frac{\frac{n}{\lambda-1}}{k}\bigg\rceil\leq \bigg\lceil \frac{t}{\lambda-1}\bigg\rceil=s<t-1~\textit{vertices~in}~V'.
\end{align}

Let $G^\star=G\cup K_r$. By \eqref{bound-r}, $G^\star$ is $d$-degenerate.
By $w_1,w_2,\ldots,w_r$, we denote the vertices of the $K_r$-component of $G^\star$.
We now extend $c'$ to a partial coloring $c^\star$ of $G^\star$ by coloring $w_i$ with $i$, where $1\leq i\leq r<k$. Clearly, $c^\star$ is an equitable tree-$k$-coloring and by \eqref{s}, 
\begin{align}\label{s-star}
  \textit{every~color~class~of}~c^\star~\textit{contains~at~most}~s+1<t~\textit{vertices~in}~V'',
\end{align}
where $V''=V'\cup \{w_1,w_2,\ldots,w_r\}$.

We color $v_{\mu+1},v_{\mu+2},\ldots,v_{n}$ one by one in such an ordering. Precisely, in the $i$-th step ($i=\mu+1,\mu+2,\ldots,n$), we attempt to color the vertex $v_{i}$ with a color from $\{1,2,\ldots,k\}$, and meanwhile, recolor (if necessary) earlier colored vertices so that the resulting coloring satisfies

(a) each color class induces a forest,

(b) each color appears on at most $t$ vertices among $V''\cup \{v_{\mu+1},\ldots,v_{i}\}$,

(c) no vertex among $V''$ is recolored.

\noindent If we succeed in a step, then turn to the next step unless all vertices of $G^\star$ have been colored, in which case the algorithm returns  TRUE, and if we fail to do such a task, the algorithm returns FALSE and ends.

In the following, we prove that the algorithm will always return TRUE. Therefore, every vertex of $G^\star$ will be colored, and by (a) and (b), the resulting coloring is an equitable tree-$k$-coloring so that each color class contains exactly $t$ vertices, since $|V(G^\star)|=kt$. Furthermore, as the vertices of the $K_r$-component of $G^\star$ have been initially colored with distinct colors and those colors would not be changed during the execution of the algorithm
by (c),
restricting the coloring of $G^\star$ to $G$ results in an equitable tree-$k$-coloring of $G$, and we win.

Suppose for contradiction that the algorithm returns FALSE at the $i$-th step with $\mu+1\leq i\leq n$. This means that $v_i$ cannot be colored (even possibly with some earlier colored vertices recolored) so that (a), (b) and (c) hold. 
The graph induced by $V''\cup \{v_{\mu+1},\ldots,v_{i-1}\}$ is denoted by $G_{i-1}$. In the following, we find contradictions by dividing the arguments into two major cases.

\textbf{Case 1. $\mu+1\leq i\leq 0.5n$.}

Let $c$ be a $k$-coloring of $G_{i-1}$ satisfying (a), (b) and (c). A color class of $c$ is \emph{big} if
it has exactly $t$ vertices, and is \emph{small} otherwise (note that there may be an empty color class). If there is a small color class of $c$ that contains at most one neighbor of $v_i$ in $G^\star$, then $v_i$ can be moved into this class, and the algorithm will be turned to the $(i+1)$-th step. Therefore, every small color class of $c$ contains at least two neighbors of $v_i$ in $G^\star$. This implies that there are $\alpha\leq \lfloor\frac{\D_{G^\star}(v_i)}{2}\rfloor=\lfloor\frac{\D_G(v_i)}{2}\rfloor$ small classes, $k-\alpha$ big classes in $c$, and at least
$t (k-\alpha)+2\alpha \geq (k-\lfloor\frac{\D_G(v_i)}{2}\rfloor)t+2\lfloor\frac{\D_G(v_i)}{2}\rfloor$ vertices have been colored under $c$. Hence by \eqref{lmdasmall} and \eqref{lmda}, we have (note that $\D_G(v_i)<\lambda d$ by the choice of $\mu$, $k\geq \frac{\Delta+1}{2}\geq \frac{9.818d+1}{2}>4.909d$, and thus $d<0.204k$)
\begin{align}
i&>\bigg(k-\bigg\lfloor\frac{\D_G(v_i)}{2}\bigg\rfloor\bigg)\cdot t+2\bigg\lfloor\frac{\D_G(v_i)}{2}\bigg\rfloor-r\\
\label{i-1} &\geq kt-(t-2)\frac{\D_G(v_i)}{2}-r \\
&>\bigg(k-\frac{\lambda d}{2}\bigg) t+(\lambda d-r)\\
\label{i-2}&>\bigg(k-\frac{\lambda d}{2}\bigg)t
\geq \left\{
\begin{array}{ll}
(k-3.5d)t>0.286kt\ge 0.286n, & \hbox{if $t\geq 4$;}\\
(k-2d)t>0.592kt\ge 0.592n, & \hbox{if $t=3$.}
\end{array}
\right.
\end{align}
Note that $\lambda d-r\geq (\lambda-\frac{1}{2})d+1>0$ by \eqref{bound-r}, \eqref{lmdasmall} and \eqref{lmda}.

From \eqref{i-2}, one immediately concludes that $i>0.592n$ if $t=3$, which contradicts the fact that $i\leq 0.5n$. Therefore, $t\geq 4$.
By \eqref{i-1},
\begin{align*}
\D_G(v_i)\cdot t>\D_G(v_i)\cdot (t-2)>2(kt-r-i)=2(n-i).
\end{align*}
On the other hand, $\D_G(v_i)\cdot t<dt(1+\frac{n}{i})$ by Lemma \ref{lem-1}. This implies
\begin{align}\label{2k-d}
  2(n-i)<dt\bigg(1+\frac{n}{i}\bigg)\leq d\bigg(\frac{n+k}{k}\bigg)\bigg(\frac{n+i}{i}\bigg).
\end{align}
Since $n> k(t-1)$ and $t\geq 4$,
\begin{align}\label{bound-k}
n+k< \frac{t}{t-1}n\leq \frac{4}{3}n.
\end{align}
Since $k> 4.909d$ and $\Delta\leq 2k-1$, by \eqref{2k-d} and \eqref{bound-k}, we conclude
\begin{align}\label{dleta-d}
  \frac{\Delta}{d}\leq \frac{2k-1}{d}<\frac{2k}{d}\leq \frac{(n+k)(n+i)}{(n-i)i}\leq \frac{4}{3}\frac{n(n+i)}{i(n-i)}=\frac{4}{3}\frac{1+\frac{i}{n}}{\frac{i}{n}\big(1-\frac{i}{n}\big)}.
\end{align}
Let $\beta=\frac{i}{n}$ and $f(\beta)=\frac{4}{3}\frac{1+\beta}{\beta(1-\beta)}$.
Since $0.286n<i<0.5n$ by \eqref{i-2} and by the condition of this case,  $0.286<\beta<0.5$.
Since $$f'(\beta)=\frac{4}{3}\frac{(1+\beta)^2-2}{\beta^2 (1-\beta)^2},$$
$f(\beta)$ decreases if $\beta<\sqrt{2}-1$ and increases if $\beta>\sqrt{2}-1$. Hence
$$f(\beta)\leq \max\big\{f(0.286),f(0.5)\big\}
<8.789,
$$
which implies $\Delta<8.789d$ by \eqref{dleta-d}, contradicting the fact that $\Delta\geq 9.818d$.

\textbf{Case 2. $0.5n< i\leq n$.}

Let $c$ be a $k$-coloring of $G_{i-1}$ satisfying (a), (b) and (c).
Define an auxiliary digraph $\mathcal{H}:=\mathcal{H}(c)$ on the color classes
of $c$ by $XY\in E(\mathcal{H})$ if and only if some vertex $x\in X\backslash V''$
has at most one neighbor in $Y$. In this case we say that $x$
\emph{witnesses} $XY$.
If $P:=M_lM_{l-1}\cdots M_1M_0$ is a path in $\mathcal{H}$ and $x_j$ ($l\geq j\geq 1$) is a vertex in $M_j$ such that $x_j$ witnesses $M_jM_{j-1}$, then switching witnesses along $P$ means moving $x_j$ to $M_{j-1}$ for every $l\geq j\geq 1$. This operation decreases
$|M_l|$ by one and increases $|M_0|$ by one, while leaving the sizes of the interior vertices of $\mathcal{H}$ (color classes) unchanged.

Let $Y_0$ be the set of small color classes of $c$, i.e, color classes having less than $t$ vertices. By $Y_j$ $(j\geq 1)$, we denote the set of color classes of $c$ such that

(i) $Y_j\cap \bigcup_{q=0}^{j-1}Y_q=\emptyset$, and

(ii) for any color class $M_j\in Y_j$ there exists a color class $M_{j-1}\in Y_{j-1}$ so that $M_{j}M_{j-1}\in \mathcal{H}$.

Let $\mathfrak{Y}=\bigcup Y_j$ and let $y=|\mathfrak{Y}|$. If $v_i$ has at most one neighbor in some color class $M_j\in Y_j\in \mathfrak{Y}$, then move $v_i$ into $M_j$ and switch  witnesses along $P=M_jM_{j-1}\cdots M_0$, where $M_t\in Y_t$ ($j\geq t\geq 0$). At this moment, $v_i$ has been colored with some earlier colored vertices recolored so that (a), (b) and (c) hold, a contradiction. Therefore, $v_i$ has at least two neighbors in every color class of $\mathfrak{Y}$, i.e.,
\begin{align}\label{a}
 \D_G(v_i)\geq 2y.
\end{align}
If there is a vertex $v_{j}\not\in \bigcup_{M\in \mathfrak{Y}}M$ with $\mu+1\leq j\leq i-1$, then
\begin{align}\label{b}
\emph{$v_j$ has at least two neighbors in every color class of $\mathfrak{Y}$,}
\end{align}
otherwise $v_j$ has at most one neighbor in some color class $M_{a}\in Y_a\in \mathfrak{Y}$ and thus $v_j\in M_{a+1}\backslash V''$, where $M_{a+1}\in Y_{a+1}\in \mathfrak{Y}$, a contradiction to the choice of $v_j$.

Divide $V(G_{i-1})\backslash \bigcup_{M\in \mathfrak{Y}}M$ into two subsets $A$ and $B$ so that $A\subset V''$ and $B\cap V''=\emptyset$. By \eqref{b}, every vertex in $B$ has at least two neighbors in every color class of $\mathfrak{Y}$, which implies
\begin{align}\label{edgepart1}
  e\bigg(B,\bigcup_{M\in \mathfrak{Y}}M\bigg)\geq 2y|B|.
\end{align}
By the choice of $V'$ and by \eqref{ordering}, $\D_G(w)\geq \lambda d$ if $w\in V'$. Therefore, the number of edges that are incident with
$A'=A\backslash\{w_1,w_2,\ldots,w_r\}$ is exactly
\begin{align}\label{edgepart2}
  \sum_{w\in A'}\D_G(w)-e(A')\geq \lambda d|A'|-e(A')> \lambda d|A'|-d|A'|=(\lambda-1)d|A'|.
\end{align}
Note that the graph induced by $A'$ is $d$-degenerate and thus $e(A')<d|A'|$.

By \eqref{a} and by Lemma \ref{lem-1}, we have
\begin{align}\label{degvi}
  2y\leq \D_G(v_i)<d(1+\frac{n}{i})<3d\leq  (\lambda-1)d
\end{align}
since $0.5n<i\leq n$ and $\lambda\geq 4$ by \eqref{lmdasmall} and \eqref{lmda}.

Note that $A\cup B=V(G_{i-1})\backslash \bigcup_{M\in \mathfrak{Y}}M$  is made up of the vertices from the $(k-y)$ color classes of $c$ that are not involved in $\mathfrak{Y}$. By the choice of $Y_0$, those $(k-y)$ color classes are big, i.e, each of them has exactly $t$ vertices. This concludes
\begin{align}\label{aandb}
|A|+|B|=t(k-y).
\end{align}
By \eqref{edgepart1}, \eqref{edgepart2} and \eqref{degvi}, we have
\begin{align}\label{edgetotal}
 e(G^\star)=e(G)+\frac{r(r-1)}{2}\geq 2y|B|+(\lambda-1)d|A'|+\frac{r(r-1)}{2}> 2y(|B|+|A|-r)+\frac{r(r-1)}{2}.
\end{align}
On the other hand, since $G$ is $d$-degenerate,  $e(G^\star)=e(G)+\frac{r(r-1)}{2}<d(kt-r)+\frac{r(r-1)}{2}$. It follows from \eqref{aandb} and \eqref{edgetotal} that
\begin{align}\label{proc}
  (2y^2-2ky+kd)t>(d-2y)r.
\end{align}
We now claim that
\begin{align}\label{yupperboundnew}
y<0.5681d.
\end{align}

Suppose first that
\begin{align*}
  \varphi(y):=2y^2-2ky+kd\geq 0.
\end{align*}
By \eqref{degvi}, we have $y< 1.5d$. Since $k\geq \frac{\Delta+1}{2}\geq \frac{9.818d+1}{2}>4.909d$,
\begin{align*}
  \varphi(1.5d)&=4.5d^2-2kd<4.5d^2-9.818d^2<0\\[.2em]
\varphi(0.5681d)&=0.64547522d^2-0.1362kd<0.64547522d^2-0.6686058d^2<0.
\end{align*}
This implies \eqref{yupperboundnew}.

On the other hand, if $\varphi(y)< 0$, then $d< 2y$ by \eqref{proc}. Using \eqref{bound-r}, \eqref{proc} and the fact that $r\geq 3$, we have
\begin{align*}
    \frac{d}{2}>r> \frac{(-2y^2+2yk-kd)t}{2y-d}\geq \frac{-6y^2+6yk-3kd}{2y-d},
\end{align*}
which implies
\begin{align*}
  \phi(y):=12y^2+(2d-12k)y+(6kd-d^2)>0.
\end{align*}
Again, by \eqref{degvi}, we have $y<1.5d$. Since $k\geq \frac{\Delta+1}{2}\geq \frac{9.818d+1}{2}>4.909d$,
\begin{align*}
  \phi(1.5d)&=29d^2-12kd<29d^2-58.908d^2<0\\[.2em]
\phi(0.5681d)&=4.00905132d^2-0.8172kd<4.00905132d^2-4.0116348d^2<0.
\end{align*}
This implies \eqref{yupperboundnew} again.

By \eqref{s} and by the condition (c) that no vertex among $V''$ would be recolored in any step, every color class of $c$ has at most $s$ vertices in $V'$. Since $A\subset V''$ and $A'\subset V'$,
\begin{align*}
  |A|=|A'|+r\leq s(k-y)+r,
\end{align*}
which implies
\begin{align}\label{blowerbound}
|B|\geq (t-s)(k-y)-r.
\end{align}
by \eqref{aandb}.

Since less than $kt$ vertices of $G^\star$ are colored and there are $k$ color classes in $c$, there is at least one color class having less than $t$ vertices, which implies that $Y_0\neq \emptyset$. Choose $M_0\in Y_0$. It follows that $|M_0|\leq t-1$.
Let $M'=M_0\backslash \{w_1,w_2,\ldots,w_r\}$ and let $l=|V'\cap M'|$. By \eqref{s}, $l\leq s$. Since $M'$ contains $|M'|-l$ vertices that are not in $V'$, each of which has degree less than $\lambda d$ by the choice of $\mu$,
the number of neighbors of $M'$ is at most
$l\Delta+(|M'|-l)\lambda d\leq s\Delta+(|M'|-s)\lambda d\leq s\Delta+(t-s-1)\lambda d$, which implies
\begin{align}\label{m0bupper}
 e(M',B)\leq s\Delta+(t-s-1)\lambda d.
\end{align}
Note that $\Delta-\lambda d\geq 9.818d-7d=2.818d>0$ by \eqref{lmdasmall} and \eqref{lmda}.

On the other hand, by \eqref{b}, every vertex in $B$ has at least two neighbors in $M_0$. Since $B\cap V''=\emptyset$, there is no edge between $B$ and $\{w_1,w_2,\ldots,w_r\}$. Hence every vertex in $B$ has at least two neighbors in $M'$, which implies
\begin{align}\label{eMB}
 e(M',B)\geq 2|B|\geq 2(t-s)(k-y)-2r>2(t-s)(k-y)-d
\end{align}
by \eqref{bound-r} and \eqref{blowerbound}. Combining \eqref{eMB} with \eqref{yupperboundnew} and \eqref{m0bupper}, we immediately have
\begin{align*}
2(t-s)\bigg(\frac{\Delta+1}{2}-0.5681d\bigg)-d\leq 2(t-s)(k-y)-d<s\Delta+(t-s-1)\lambda d,
\end{align*}
which implies
\begin{align}\label{basic}
(t-2s)\Delta<\bigg(1.1362(t-s)+(t-s-1)\lambda+1\bigg)d.
\end{align}

We now divide the proofs into two subcases according to the value of $t$.

\textbf{Subcase 2.1. $3\leq t\leq 6$.}

In this case $s=1$ by \eqref{s(t)} and $\lambda=1+t$ by \eqref{lambda}. Since $\Delta\geq 9.818d$, we can deduce from \eqref{basic} that
\begin{align*}
9.818(t-2)<1.1362(t-1)+(t-2)(t+1)+1,
\end{align*}
which implies $t\geq 8$, a contradiction.

\textbf{Subcase 2.2. $t\geq 7$.}

In this case $\frac{t}{s}=\lambda-1$ by \eqref{lambda}.
 Since $\Delta\geq 9.818d$, we can deduce from \eqref{basic} that
\begin{align*}
9.818\leq \frac{\Delta}{d}<\frac{t-s}{t-2s}(1.1362+\lambda)=
\frac{\frac{t}{s}-1}{\frac{t}{s}-2}(1.1362+\lambda)=\frac{(\lambda-2)(\lambda+1.1362)}{\lambda-3}:=\theta(\lambda).
\end{align*}
Since $\theta(4.1818)<9.81792$ and $\theta(6)<9.515$, 
$\lambda<4.1818$ or $\lambda>6$, contradicting \eqref{lmda}. \hfill$\square$

\section{Remarks and Problems for Further Research}

We note that by being a little more careful the condition $k\geq (\Delta+1)/2$ in Theorem \ref{main} can be replaced by $k\geq \Delta/2$.
With this result in mind, we conjecture that the bound in Conjecture \ref{conj1} can be improved for most connected graphs with
maximum degree that is even.

\begin{conj}\label{conjnew}
Every graph $G$ is equitably tree-$\lceil\Delta(G)/2\rceil$-colorable unless  $G$ is a cycle or a complete graph on odd number of vertices.
\end{conj}

As referred in Section 1, Esperet, Lemoine and Maffray \cite{ELM.2015} showed that $va_{eq}^*(G)\leq 3^{d-1}$ for every $d$-degenerate graph $G$. Actually, it is interesting to find an $O(d)$ upper bound  (partial result can be found in \cite{ZN19}). Recently, Li and Zhang \cite{LZ20} proved that $va_{eq}^*(G)\leq d$ for every graph $G$ with treewidth at most $d$, and showed that the upper bound $d$ of the equitable vertex arboreal threshold is sharp. Indeed, there are graphs with treewidth $d$ that are not equitably tree-$(d-1)$-colorable, see \cite{LZ20}.
Since every graph with treewidth at most $d$ is $d$-degenerate, we naturally release the following conjecture.

\begin{conj}\label{conj-deg}
Every $d$-degenerate graph $G$ is equitably tree-$d'$-colorable for any integer $d'\geq d$, i.e, $va_{eq}^*(G)\leq d$.
\end{conj}

Another interesting problem is to lower the bound $9.818$ of $\Delta/d$ in Theorem \ref{main}. If it can be improved to $1$, then Equitable Vertex Arboricity Conjecture (Conjecture \ref{conj1}) will be verified. This can be easily seen from a trivial fact that every graph with maximum degree at most $\Delta$ is $\Delta$-degenerate.
For Conjecture \ref{conj3}, it may be also interesting to investigate the following problem firstly, instead of proving Conjecture \ref{conj3} directly.

\begin{problem}
Find a constant $\beta\leq 1$ (the larger the better) such that every $n$-vertex planar graph with maximum degree at most $\Delta$ is equitably tree-$3$-colorable provided that $\Delta\leq \beta n$.
\end{problem}

\noindent Clearly, if $\beta=1$ is an affirmative answer to the above problem, then Conjecture \ref{conj3} will be verified.

No suprisely, the list version of Equitable Vertex Arboricity Conjecture \cite{KMP20,Z16,DDFS18,DFS20,DFS20-conf,LZ2020} and the list version of all problems considered or mentioned in this paper are also interesting.

\bibliographystyle{srtnumbered}
\bibliography{mybib}

\end{document}